\documentclass[12pt]{amsart}
\usepackage{amsthm}
\usepackage{amssymb}
\usepackage{amsmath}
\theoremstyle{plain}
\newtheorem{proposition}{Proposition}
\newtheorem{theorem}[proposition]{Theorem}
\newtheorem{lemma}[proposition]{Lemma}
\newtheorem{corollary}[proposition]{Corollary}

\theoremstyle{definition}
\newtheorem{definition}[proposition]{Definition}

\theoremstyle{definition}

\newtheorem{remark}[proposition]{Remark}

\numberwithin{equation}{section}

\numberwithin{proposition}{section}

\gdef\myletter{}

\let\savetheequation\theequation
\def\theequation{\savetheequation\myletter}

\newcommand{\CC}{{\mathbb C}}

\newcommand{\RR}{{\mathbb R}}

\renewcommand{\date}{\today}
\def \bar{\overline}

\begin{document}

\author{Thomas Bloom*}
\thanks{*supported by an NSERC of Canada grant}
\address{University of Toronto, Toronto, Ontario, M5S2E4, Canada}
\email{bloom@math.toronto.edu}
\title{Almost sure convergence for Angelesco ensembles}

\maketitle 

\vspace{-4mm}

\begin{center}
June 20, 2012
\end{center}

\vspace{4mm}

\begin{abstract}
Almost surely in an Angelesco ensemble, the normalized counting measure of a random point converges weak* to the equilibrium measure.We also establish a large deviation principle with good rate function and unique minimizer.
A large deviation principle for Angelesco ensembles was established by Eichelsburger,Sommerauer, and Stolz \cite{ESS},however,the method of proof given here is different and settles a question of Kuijlaars\cite{K}
in full generality. 
\end{abstract}
\section{introduction}
In this paper we will establish almost sure convergence results for the normalized counting measure of a random point in an  Angelesco ensemble.We will also establish a large deviation principle (l.d.p.) for such ensembles.  Angelesco ensembles are special cases of multiple orthogonal polynomial ensembles.\\

First we review the corresponding theory for what, in this context, should be called orthogonal polynomial ensembles.
 The Gaussian Unitary ensemble has the joint probability distribution (j.p.d.)\\
\begin{equation}Prob_n=\frac{A_n}{Z_n}\mathrm{d}x_1\ldots\mathrm{d}x_n\end{equation}\\where
\begin{equation}A_n=\prod_{1\leq i<j\leq n} |x_i-x_j|^2 \exp (-n\sum_{i=1}^n x_i^2)\end{equation}where $Z_n$ is a constant (which is explicitly known in this case). It
is known that the random empirical measure of the eigenvalues (i.e., the normalized counting measure of a random point under the above probability measure) converges, almost
surely to the (non-random) measure (known as the Wigner semi-circle law)
$$\mu=\frac{1}{\pi}\sqrt{2-x^2}$$ for $-\sqrt{2}\leq x\leq \sqrt{2}$.\\

The above measure may be interpreted in terms of potential theory. It is the equilibrium measure of $\RR \subset \CC$ with weight $\exp(\frac{-x^2}{2})$. That is $\mu$ is the unique probability measure on $\RR$ which minimizes the functional\begin{equation}E(\nu)=\int_{\RR}\int_{\RR}\log\frac{1}{|x-y|} \mathrm{d}\nu (x) \mathrm{d}\nu (x)+\int_{\RR}x^2\mathrm{d}\nu \end{equation}over probability measures on $\RR$. $\mathrm{E}(\nu)$ is referred to as the weighted energy of $\nu$\\

A general unitary ensemble has the joint eigenvalue probability distribution
 on $ \RR^n$ given by\begin{equation}Prob_n=\frac{A_n}{Z_n}\mathrm{d}x_1\ldots\mathrm{d}x_n\end{equation}\\where\begin{equation}A_n = 
\prod_{1\leq i<j\leq n}|x_i -x_j|^2 \exp (-2n(\sum_{i=1}^n Q(x_i))) \end{equation}where $Z_n$ is a constant and $Q$ is a continuous function on $\RR$ with appropriate growth at $\infty$. Again the random empirical measure of the eigenvalues converges, almost surely to a non-random measure $\mu_Q$, which may be characterized in terms of potential theory -it is the weighted equilibrium measure of $\RR \subset \CC$ with weight $\exp (-Q)$ and, as such minimizes the functional\begin{equation}E^Q(\nu)=\int_{\RR}\int_{\RR}\log\frac{1}{|x-y|} \mathrm{d}\nu (x) \mathrm{d}\nu (x)+2\int_ {\RR}Q(x)\mathrm{d}\nu \end{equation}over probability measures on $\RR$.\\

A large deviation principle (l.d.p.) is also known to hold for general Unitary ensembles. The speed of the l.d.p. is $n^2$, there is a good rate function $E^Q(\mu)-E^Q(\mu_Q)$ and $\mu_Q$ is the unique minimizer of $E^Q$ and so also of the rate function. The l.d.p. for the Gaussian Unitary ensemble was established by Ben Arous and Guionnet(\cite{BG}) and relies on earlier work of Voiculescu. It is well-known (see \cite{HP}, page 211) that the almost sure convergence of the normalized counting measure of a random point follows from 
the l.d.p. An l.d.p. has also been established for numerous other ensembles ($\beta$-ensembles, ensembles in the complex plane,etc..-see \cite{HP} or\cite{AGZ}) by various authors using extensions of the method of Ben Arous-Guionnet. Furthermore, recently an   l.d.p. for Angelesco (and other) ensembles was proved by Eichelsbacher, Sommerauer and Stolz \cite{ESS} again extending the Ben Arous-Guionnet method.\\

In recent years, using new developments in pluripotential theory due primarily to R.Berman and S.Boucksom a l.d.p.(see \cite{Be} and the references given there) has been established  for certain canonical point processes on complex manifolds (which may be viewed as a multivariable generalization of the case of Unitary ensembles). The methods developed there (see also\cite{BL}) give  new results and a different approach even in the one variable case (i.e. Unitary ensembles). In this paper we adapt and simplify these methods to the one-variable case, in particular Angelesco ensembles. We remark that in the case of a single interval, $\Gamma$, an Angelesco ensemble reduces to a Unitary ensemble on that interval, so, in fact, Unitary ensembles can be considered as a special case of Angelesco ensembles and the methods of this paper appply to Unitary ensembles.

Angelesco ensembles arise in approximation problems\cite{A}
and also some aspects of Brownian motion\cite{DK}.\\

In this paper, even though the problems can be entirely formulated in terms of real variables, we use potential theory in the complex plane. In fact it is essentially exclusively potential theory that we use.\\

 We show (theorem 4.4) that almost surely in an Angelesco ensemble, the normalized counting measure of a random point converges weak* to the equilibrium measure. Then in theorem 5.6 we show that the same result holds for weighted Angelesco ensembles. As in the case of Unitary ensembles, the limiting measure is non-random.\\

The problem of providing a rigorous proof of this result was raised by Kuijlaars\cite{K}.The
 method of \cite{ESS} requires the restriction that the functions $w_i$, occuring in the definition of an Angelesco ensemble (see definition 2.1) be continuous on $\Gamma_i$. In this paper the results are established only requiring that the measures $w_i dx$ satisfy the Bernstein-Markov (BM) inequality (see definitions 3.8 and 5.2). This class of measures is sufficiently general to include the case $w_i>0$ a.e on the interval $\Gamma_i$ and so gives a rigorous proof of the almost sure convergence in the generality posed by Kuijlaars.\\

As mentioned, the method of proof we use is different from that of Ben Arous-Guionnet even when applied to Unitary ensembles. First we define and establish properties of Fekete points in an Angelesco ensemble (corollaries 3.6,3.7). Then we prove a type of large deviation estimate (theorem 4.1) where the measures $w_i(x)dx$ need only satisfy a BM inequality  (this  large deviation estimate is weaker than the full l.d.p.-it is a type of estimate whose utility was shown by Johannson \cite{J} and we refer to it as a Johannson large deviation). This is then used to prove theorem 4.4 which establishes the almost sure convergence of the normalized counting measure of a random point in an Angelesco ensemble. In section 5 the same method is used in the case of weighted Angelesco ensembles and the measures $w_i\mathrm{d}x$ must satisfy a weighted BM inequality(see definition 5.2).

  We state a (full) l.d.p. for Angelesco ensembles (theorem 6.2). The proof is given in section 7. Again the method of proof is different (even in the case of unitary ensembles) from the method of Ben Arous-Guionnet and it is valid for measures $\tau_i$ satisfying a strong BM inequality(see definition 5.3  ). The collection of such measures is again known to include the case $w_i\mathrm{d}x$ where $w_i >0$ a.e. on an interval.\\To prove the (full) l.d.p.
first we use the Johannson large deviation result to prove the l.d.p. result for measures which are equilibrium measures (corollary 7.5). Then the result
is extended to the general case by appropriately approximating a general measure by an sequence of equilibrium measures (theorem 7.6).\\

$\mathbf{Acknowledgement}$.The research for this paper originated at the conference on New Perspectives in Univariate and Multivariate Orthogonal Polynomials, held at Banff International Research Station, October 10-15, 2010. I would like to thank BIRS for the excellent support. I would also like to thank A.Kuijlaars and N.Levenberg  for helpful conversations.

\section{Angelesco ensembles}
In this section we will describe Angelesco ensembles. They are, in fact, a particular case of multiple orthogonal polynomial (MOP) ensembles  which we will outline first. We will closely follow the presentation of Kuijlaars\cite{K}.
For more details and extensive references on MOP's see \cite{K} or \cite{K1}.
\begin{definition}Given integrable a.e. positive functions $w_1,\ldots,w_p$ on $\RR$ and a multi-index $\vec{n}=(n_1,\ldots,n_p) \in \mathbb{N}^p$, a MOP is a monic polynomial $P_{\vec{n}}$
of degree $n=|\vec{n}|=n_1+\ldots +n_p$, such that \begin{equation}\int_{-\infty}^{+\infty}
P_{\vec{n}}(x)x^kw_j(x)\mathrm{d}x=0\end{equation}
 for $k=0.\ldots,n_j-1;j=1,\ldots,p$.\end{definition}(the integrals are assumed to be absolutely convergent)

In \cite{K} these are referred to as type II MOP's, but, this being the only type used in this paper we will simply refer to them as MOP's. Also in that paper the functions $w_i$ are referred to as weight functions, but we will reserve that term for what in \cite{K} are called "varying weights".

The above system gives $n$ linear equations for the $n$ free coefficients of the monic polynomial $P_{\vec{n}}$. The equations do not always have a unique solution, but in the Angelesco case they do. The Angelesco case requires the following data:\\ - $p$ disjoint (non-degenerate) compact intervals $\Gamma_1,\ldots,\Gamma_p\subset\RR$ \\(We assume that $\Gamma_{j+1}$ lies to the right of $\Gamma_j$ for $j=1,\ldots,(p-1)$.) \\ - integrable a.e. positive functions $w_1,\ldots,w_p$ satisfying $supp(w_i)\subset\Gamma_i$ for $i=1,\ldots,p$. That is $w_i\equiv 0$ on 
$\RR\setminus\Gamma_i$.\\

 Associated to this situation there is a probability measure on $\RR ^n$, denoted $Prob_
{\vec{n}}$ and characterized as follows:\\

Given $X=(x_1,\ldots,x_n)\in\RR^n$ we let $\rho X$ denote the unique element of $L_n$ obtained from $X$ by permutation of coordinates,where\\ \begin{equation}L_n=\{(x_1,\ldots,x_n)\in\RR^n |x_1\leq x_2\ldots\leq x_n\}\end{equation}  So\begin{equation}\rho:\RR^n\to L_n.\end{equation}

We use the notation$$\Delta(X)=\prod_{1\leq j<k\leq n}(x_k-x_j)$$ where $X=(x_1,\ldots,x_n)$ and $$\Delta(X,Y)=\prod_{k=1}^n\prod_{j=1}^m(x_k-y_j)$$ for $X=(x_1,\ldots,x_n)$ and $Y=(y_1,\ldots,y_m)$.\\
$Prob_{\vec{n}}$ is invariant under permutation of coordinates,
so it is characterized by its push-forward by $\rho$ given below.\\

We set $N_j=\sum_{i=1}^j n_i$ and $N_0=0$.\\Also we set$$x_k^{(j)}=x_{N_{j-1}+k}$$ for $k=1,\ldots,n_j;j=1,\ldots ,p.$ Then,\begin{equation}\rho_*(Prob_{\vec{n}})=\frac{A_{\vec{n}}}{Z_{\vec{n}}}\prod_{i=1}^p\prod_{k=1}^{n_i}w_i(x_k^{(i)})\mathrm{d}x_k^{(i)}.\end{equation}Note that the support of $Prob_{\vec{n}}$ is $\{X\in\RR^n|\rho(X)\in\Gamma_1^{n_1}\times\ldots\times\Gamma_p^{n_p}\subset L_n\}$\\
Here $Z_{\vec{n}}$ is a normalizing constant \begin{equation}Z_{\vec{n}}=\int_{L_n}A_{\vec{n}}(X)\prod_{i=1}^p\prod_{k=1}^{n_i}w_i(x_k^{(i)})\mathrm{d}x_k^{(i)},\end{equation}and \begin{equation}A_{\vec{n}}(X)=\prod_i^p\Delta (X^{(i)})^2\prod_{1\leq i<j\leq p}\Delta(X^{(j)},X^{(i)}),\end{equation} where $X^{(i)}=(x_1^{(i)},\ldots,x_{n_i}^{(i)})$.\\

 Note that $A_{\vec{n}}$ is $\geq 0$, by our convention on labelling the $\Gamma_j.$\\ 

For the above probability distribution  $\mathbb{E}(\prod_{j=1}^n (z-x_j))=P_{\vec{n}}(z).$ That is the monic polynomial $P_{\vec{n}}(z)$ is the expected value of the polynomial whose roots are a random point in the Angelesco ensemble \cite{K}.\\

A point $X$ in an Angelesco ensemble (of order $\vec{n}$) is a point of $\RR^n$ in the support of $Prob_{\vec{n}}$ (Hence in $\Gamma_1^{n_1}\times\ldots\times\Gamma_p^{n_p}$ or an image of this set under a permutation of the coordinates of $\RR^n$).\\
 
The normalized counting measure of the point $X$, denoted $\frac{1}{n}\delta (X)$, is $\frac{1}{n}\sum_{i=1}^n\delta (x_i)$ ,($\delta$ denoting the Dirac measure) where $X=(x_1,\ldots ,x_n)$.  This is a discrete measure with support in $\Gamma_1\cup \ldots \cup\Gamma_p$. It is independent of any permutation of the coordinates of $X$.
  That is $$ \frac{1}{n}\delta(X)=\frac{1}{n}\delta(\rho X).$$ 

 Given positive reals $r_1,\ldots,r_p$ with \begin{equation}\sum_{i=1}^p r_i=1\end{equation} we consider p-tuples of positive Borel measures $(\nu_1,\ldots,\nu_p)$ with $\nu_i$ supported on $\Gamma_i$ and \begin{equation}\int_{\Gamma_i}\mathrm{d}\nu_i=r_i\end{equation} for $i=1,\ldots,p$.\\
Associated to the data of an Angelesco ensemble, and $r=(r_1,\ldots,r_p)$ as in (2.7) and there is an energy minimizing problem.
Let $\mathcal{M}(\Gamma)$ denote the space of probability measures on $\Gamma=\Gamma_1\cup\ldots\cup\Gamma_p$ and let $\mathcal{M}_r(\Gamma)$ denote the p-tuples of positive Borel measures $\nu=(\nu_1,\ldots,\nu_p)$ with $supp(\nu_i) \subset \Gamma_i$ for $i=1,\ldots,p$ and satisfying (2.7) and (2.8). Then  $\mathcal{M}_r(\Gamma)\subset\mathcal{M}(\Gamma)$ and we consider these spaces to have the weak* topology. By a (particular case) of a theorem of Gonchar-Rachmanov \cite{GR}, there is a unique element of $\mathcal{M}_r(\Gamma)$, denoted $\mu_{\Gamma}=(\mu_{\Gamma}^1,\ldots,\mu_{\Gamma}^p)$ which minimizes the following functional $E$ ( we will refer $\mu_{\Gamma}$ to as the equilibrium measure). $E$ is defined as follows: \begin{equation}E(\nu)=E(\nu_1,\ldots,\nu_p)\end{equation}\\$$=\sum_{i=1}^p \int_{\Gamma_i}\int_{\Gamma_i}\log \frac{1}{|x-y|}\mathrm{d}\nu_i(x)\mathrm{d}\nu_i(y)+\sum_{1\leq i<j\leq p}\int_{\Gamma_i}\int_{\Gamma_j}\log\frac{1}{|x-y|}\mathrm{d}\nu_i(x)\mathrm{d}\nu_j(y).$$That is $E(\mu_{\Gamma})=\inf_{\nu\in\mathcal{M}_r(\Gamma)}E(\nu)$.\\

Of course  $\mu_{\Gamma}$ depends on $r_1,\ldots,r_p$ although our notation does not explicitly indicate this.

\section{Fekete points}

The object of this section is to define Fekete points for Angelesco ensembles and to establish their asymptotic properties. In the case of unitary ensembles, Fekete points  (referred to as weighted Fekete points) of order $n$ are the coordinates of a point in $\RR^n$ at which $A_n$ (see (1.5)) assumes its maximum on $\RR^n$. We define Fekete points analogously for the Angelesco ensemble. Specifically we have
\begin{definition}Let $\vec{n}=(n_1,\ldots,n_p)$ be a p-tuple of positive integers. Let $F\in supp(Prob_{\vec{n}}),$
 We say that $F$ is a Fekete point (for an Angelesco ensemble) of order  $\vec{n}=(n_1,\ldots,n_p)$ if \begin{equation}A_{\vec{n}}\rho (F)
\geq A_{\vec{n}}(X) \end{equation}for all $ X\in\Gamma_1^{n_1}\times\ldots\times\Gamma_p^{n_p}$ where $A_{\vec{n}}$ is given by (2.6). 
 \end{definition}
As with unitary ensembles Fekete points of order $\vec{n}$ are not necessarily unique but nevertheless, the normalized counting measures of a sequence of  Fekete points converge as follows:

We consider  sequences of multi-indices  $\vec{n}(d)$, for $d=1,\ldots$ for which\begin{equation}\lim_{d \to \infty}\frac{\vec{n}(d)}{n(d)}=(r_1,\ldots,r_p) \end{equation} where each $r_i>0.$ We let $F(\vec{n}(d))$ denote a corresponding sequence of Fekete points. The theorem below shows that the  discrete measures formed by the normalized counting measures of the Fekete points converge, weak* to the equilibrium measure.
\begin{theorem}$\lim_d\frac{1}{n(d)}\delta(F(\vec{n}(d))=\mu_{\Gamma}$, weak*.\end{theorem}
We begin with two lemmas.
\begin{lemma}The functional $-E$ (given by (2.7))is 
uppersemicontinuous (usc) on $\mathcal{M}.$\end{lemma}
\begin{proof}$\log|x-y|$ is usc so by a basic property of usc functions (\cite{ST},\\chapter0, theorem I4) it follows that $-E$ is usc.\end{proof}
\begin{lemma}Let $X_{\vec{n}(d)}$ for $d=1,2,\ldots$ be a sequence of points in an Angelesco ensemble with $\vec{n}(d)$ satisfying (3.2). Suppose that the sequence of normalized counting measures of the points converges to $\sigma \in \mathcal{M}_r(\Gamma).$ Then$$\limsup \frac{1}{(n(d))^2}\log A_{\vec{n}(d)}(\rho X_{\vec{n}(d)})\leq -E(\sigma).$$\end{lemma}
\begin{proof}Note that if $\frac{1}{n}\sum_{j=1}^n\delta(p_j)$ is a sequence of discrete measures on a compact set $K\in\RR^n$ converging, weak*, to a measure $\tau$, then $\frac{1}{n^2}\sum_{j,k=1,j\ne k}^{n}\delta (p_j,p_k)$ converges to the product measure $\tau \otimes \tau$ on $K\times K$.\\  Also$$\log A_{\vec{n}}(X)=\sum_{i=1}^p\sum_{\substack{r,s=1\\r\neq s}}^{n_i}\log|x_r^{(i)}-x_s^{(i)}|+\sum_{1\leq i<j\leq p}\sum_{r=1}^{n_i}\sum_{s=1}^{n_j}\log|x_r^{(i)}-x_s^{(j)}|.$$ so that $-\log A_{\vec{n}}$ is a discrete version of the energy given by (2.5). The result now follows from lemma 3.3.
\end{proof}
\begin{corollary}Let $X_{\vec{n}(d)}$ be a sequence of points in an Angelesco ensemble with $\vec{n}(d)$ satisfying
(3.2). Then $$\limsup_d\frac{1}{(n(d))^2}\log A_{\vec{n}(d)}(\rho X_{\vec{n}(d)})\leq -E(\mu_{\Gamma}).$$\end{corollary}

\begin{proof}(of theorem 3.2)
Consider the sequence of discrete measures $\frac{1}{n(d)}\delta (F(\vec{n}(d))$.

Since the weak* topology is metrizable and $\mathcal{M}(\Gamma)$ is compact, it suffices to show that any subsequence of $\frac{1}{n(d)}\sum\delta (F(\vec{n}(d)))$ has a further subsequence which converges to $\mu_{\Gamma}$. Thus we may assume that the sequence of of discrete measures converges to $\sigma=(\sigma_1,\ldots,\sigma_p) \in \mathcal{M}_r(\Gamma).$ Since the equilibrium measure is unique and characterized by minimizing the functional $E$, to show that $\sigma =\mu_{\Gamma}$ it suffices to show that $E(\sigma )=E(\mu_{\Gamma} ).$
By corollary 3.5
 it remains to show that $$\liminf_d\frac{1}{(n(d))^2}\log A_{\vec{n}(d)}\rho (F(\vec{n}(d)) \geq -E(\mu_{\Gamma}).$$

Consider probability measures $\lambda_i$ on $\Gamma_i$, for $i=1,\ldots,p$. It will suffice to show that $$\liminf_d\frac{1}{(n(d))^2}\log A_{\vec{n}(d)}\rho (F(\vec{n}(d)))\geq -E(\bar{\lambda}_1,\ldots,\bar{\lambda}_p)$$ where $\bar{\lambda}_i=r_i \lambda_i$, since an arbitrary element of $\mathcal{M}_r(\Gamma)$ is of the form $(r_1\lambda_1,\ldots,r_p\lambda_p)$.\\

 After taking logarithms in (2.6) and integrating over $\Gamma_1^{n_1}
\times \ldots \times\Gamma_p^{n_p}$ with respect to $ \prod_{i=1}^p\prod_{j=1}^{n_i} \mathrm{d}\lambda_j(x_j^{(i)})$ we obtain, from the defining property of Fekete points \begin{equation}\log A_{\vec{n}(d)}\rho (F(\vec{n}(d)))\geq \sum_{i=1}^p\sum_{\substack{r,s=1\\r\ne s}}^{n_i}\int_{\Gamma_i}\int_{\Gamma_i}\log|x_r^{(i)}-x_s^{(i)}|\mathrm{d}\lambda_i(x_r^{(i)})\mathrm{d}\lambda_i(x_s^{(i)}) \end{equation} $$+ \sum_{1\leq i<j\leq p}\sum_{r=1}^{n_i}\sum_{s=1}^{n_j}\int_{\Gamma_i}\int_{\Gamma_j}\log|x_r^{(i)}-x_s^{(j)}| \mathrm{d}\lambda_i(x_r^{(i)}) \mathrm{d}\lambda_j(x_s^{(j)})$$=\begin{equation}\sum_{i=1}^p n_i(n_i -1)\int_{\Gamma_i}\int_{\Gamma_i} \log|x-y|\mathrm{d}\lambda_i(x)\mathrm{d}\lambda_i(y)\end{equation}$$ + \sum_{1\leq i<j\leq p}n_i n_j\int_{\Gamma_i}\int_{\Gamma_j}\log|x-y|\mathrm{d} \lambda_i(x)\mathrm{d}\lambda_j(y)$$
We may assume that $E(\bar{\lambda})<+\infty$ (since if $E(\bar{\lambda})=+\infty$,the inequality we wish to prove certainly holds) so that all the integrals for $i,j=1,\ldots,p.$\\$$ \int_{\Gamma_i}\int_{\Gamma_j}\log|x-y|\mathrm{d} \lambda_i(x)\mathrm{d}\lambda_j(y)$$ are bounded. Hence using the fact that (3.2) is satisfied, 
the right hand side of (3.4) is, after dividing by $(n(d))^2$, for any $\epsilon > 0$ and $d$ sufficiently large$$\geq -E(\bar{\lambda})-\epsilon.$$ $\epsilon$ being arbitrary, the result is proved.\end{proof}
\begin{corollary}$$\lim_d\frac{1}{(n(d))^2}\log A_{\vec{n}(d)}\rho F(\vec{n}(d))=-E(\mu_{\Gamma})$$\end{corollary}
Let $X_{\vec{n}(d)}$ be a sequence of points in an Angelesco ensemble (we always assume that the sequence of multi-indices $\vec{n}(d)$ satisfies (3.2)). We say that the sequence is $asymptotically\, Fekete$ if $$\lim_d\frac{1}{(n(d))^2}\log A_{\vec{n}(d)}(\rho X_{\vec{n}(d)})=-E(\mu_{\Gamma})$$.\begin{corollary}Let $X_{\vec{n}(d)}$ be an asymptotically
Fekete sequence of points in an Angelesco ensemble.Then $\frac{1}{n(d)}\sum\delta(X_{\vec{n}(d)})$ converges 
to $\mu_{\Gamma}$ weak*.\end{corollary}\begin{proof}Passing to a subsequence, if necessary, we may assume that $\frac{1}{n(d)}\sum\delta(X_{\vec{n}(d)})$ converges to $\sigma\in\mathcal{M}_r(\Gamma)$. Then \begin{equation}\limsup_d\frac{1}{(n(d))^2}\log A_{\vec{n}(d)}(\rho X_{\vec{n}(d)})\leq -E(\sigma)\leq -E(\mu_{\Gamma}).\end{equation}But, by  hypothesis, $X_{\vec{n}(d)}$ is asymptotically Fekete, so we must have equalities in (3.5)  above and therefore $\sigma =\mu_{\Gamma}$\end{proof}We will also
need to establish the asymptotic behaviour of the normalizing constants $Z_{\vec{n}(d)}$.We will use the concept of a positive measure satisfying a Bernstein-Markov (BM) inequality, defined below.\\
\begin{definition}
A measure $\tau$ on a compact subset $K$ of $\CC$ satisfies the \emph{BM inequality}, if given $\epsilon
 > 0$ there exists a constant $C >0$ such that $$\sup_K|p| \leq C(1+\epsilon)^{deg(p)}\int_K|p|\mathrm{d}\tau$$ for all analytic polynomials $p$.\end{definition}
The terminology "regular measure" is also used for this property\cite{StT}.  It is also common to define the BM inequality with the $L^2$ norm on the right side of the definition, but, in fact,  it is equivalent to use any $L^p$ norm in the definition.  (\cite{StT},theorem 3.4.3).\\

It is known that if $w> 0$ a.e. on an interval,  the measure $w(x)\mathrm{d}x$ satisfies the BM property on that interval. (see\cite{StT}, corollary  4.1.3). It is also known that there is a large class of measures which satisfy the BM inequality on an interval \cite{StT}. In fact, there are measures, not absolutely continuous with respect to Lebesgue measure
but which still satisfy the BM inequality\cite{BL},\cite{StT}\\

Henceforth we will deal with general measures satisfying the BM inequality. That is\begin{equation}\rho_*(Prob_{\vec{n}})=
\frac{A_{\vec{n}}(X)}{Z_{\vec{n}}}\mathrm{d}\tau (X) \end{equation}where\begin{equation}\mathrm{d}\tau (X)=\prod_{i=1}^p\prod_{k=1}^{n_i}\mathrm{d}\tau_i(x_k^{(i)})\end{equation}and each $\tau_i$ is a measure on $\Gamma_i$
satisfying the BM inequality on that interval  and also (2.7),(2.8) are satisfied.\\
Furthermore \begin{equation}Z_{\vec{n}}=\int_{L_n} A_{\vec{n}}(X)\mathrm{d}\tau (X).\end{equation}

 $A_{\vec{n}(d)}$ is a polynomial of degree $\leq 2n$ in each variable (actually it is of degree $n_i +n$ in $x_k^{(i)}$). Each measure $\mathrm{d}\tau_i$ satisfies the BM inequality on $\Gamma_i$ so we can, given $\epsilon$, choose the constant $C$ in  definition 3.8 valid for all $i=1,\ldots,p$\begin{theorem}$\lim_d\frac{1}{(n(d))^2}\log Z_{\vec{n}(d)}=-E(\mu_{\Gamma})$\end{theorem}\begin{proof}First we note that by corollary 3.6, the result is true with $A_{\vec{n}(d)}\rho F(\vec{n}(d))$ replacing $Z_{\vec{n}(d)}$ so we need only show that these two quantities have the same asymptotic behaviour.\\
Let$$\rho(F)=(f_j^{(i)}),j=1,\ldots,n_i;i=1,\ldots,p$$
Now,$A_{\vec{n}(d)}\rho (F) = \sup_{x_1^{(1)} \in \Gamma^1}A_{\vec{n}(d)}(x_1^{(1)},f_2^{(1)},\ldots,f_{n_p}^{(p)})$ which  is, by the BM inequality $$\leq 
C(1+\epsilon)^{2n}\int_{\Gamma_1}A_{\vec{n}(d)}(x_1^{(1)},f_2^{(1)},\ldots,f_{n_p}^{(p)})\mathrm{d}\tau_1(x_1^{(1)})$$By repeating the use of the BM inequality in each variable,we obtain$$A_{\vec{n}(d)}(F(\vec{n}(d))\leq C^n(1+\epsilon)^{2n^2}\int_{\Gamma_1^{n_1}}\ldots \int_{\Gamma_p^{n_p}}A_{\vec{n}(d)}\mathrm{d}\tau (X)$$ $$=C^n(1+\epsilon)^{2n^2}Z_{\vec{n}(d)}.$$Since $\epsilon >0$ is arbitrary and (3.2) is satisfied,we have $$-E(\mu_{\Gamma})=\lim_d
\frac{1}{(n(d))^2}\log A_{\vec{n}(d)}\rho (F(\vec{n}(d))\leq \limsup_d\frac{1}{(n(d))^2}\log Z_{\vec{n}(d)}.$$For the inequality in the other direction,simply note that it follows from the defining equation for $Z_{\vec{n}}$ (see (2.3)), that we get $Z_{\vec{n}} \leq A_{\vec{n}}\rho (F(\vec{n}))C^n$ where $C\geq \int_{\Gamma_i}\mathrm{d}\tau_i$ for $i=1,\ldots,p$.\end{proof}
\section{a.s. convergence}
In this section we will establish almost sure weak* convergence of the normalized counting measure of a random point in an Angelesco ensemble. First  we will use a type of large deviation estimate, weaker than a full l.d.p., used by Johannson\cite{J} which when combined with results of section 3 shows that almost surely a sequence of points in an Angelesco ensemble is
asymptotically Fekete. By the (deterministic) result corollary 3.7 the normalized counting measure of asymptotically Fekete sequences converges weak* to the equilibrium measure and Theorem 4.4 follows. \\

We will begin with a version of the Johannson large deviation\cite{J}. On $\RR^n$, let $\lambda_n$ be a finite measure with compact support, $K_n$,  let $g_n$ be a non-negative, continuous function and$$Z_n=\int_{K_n}g_n\mathrm{d}\lambda_n$$for $n=1,2,\ldots$such that
$$\lim_{n \to \infty}\frac{1}{n^2}\log\int_{K_n}g_n\mathrm{d}\lambda_n(x)=\lim_{n\to\infty}\frac{1}{n^2}\log\sup_{K_n}g_n=\log\gamma >-\infty$$\\for some constant $\gamma.$
Given $\eta >0$ let $$B_{\eta ,n}=\{x^{(n)} \in K_n|g_n^{\frac{1}{n^2}}\leq \gamma -\eta\}.$$

Then we have:\begin{theorem}$$\frac{1}{Z_n}\int_{B_{\eta ,n}}g_n\mathrm{d}\lambda_n\leq (1-\frac{\eta}{2\gamma})^{n^2}$$ for $n$ sufficiently large.\end{theorem}
\begin{proof}By hypothesis,$$\lim_nZ_n^{\frac{1}{n^2}}=\gamma.$$So given $\epsilon >0$, we have$$Z_n\geq (\gamma -\epsilon)^{n^2}$$ for $n$ sufficiently large. The defining property of the set $B_{\eta ,n}$ yields $$\frac{1}{Z_n}\int_{B_{\eta ,n}}g_n\mathrm{d}\lambda_n\leq \frac{(\gamma -\eta)^{n^2}}{(\gamma -\epsilon)^{n^2}}$$ for $n$ sufficiently large.\\
But $\frac{\gamma -\eta}{\gamma -\epsilon}\leq 1-\frac{\eta}{2\gamma}$ for $\epsilon$ sufficiently small (given $\eta$).The result follows.\end{proof}
We consider the probability spaces $$(K_n,\frac{g_n\lambda_n}{Z_n}).$$ \\Then theorem 4.1 says that given $\eta >0$ then the probability of $B_{\eta ,n}$ is $\leq\exp(-cn^2)$ for some constant $c>0$ and $n$ sufficiently large.\\Next we consider the product probability space $$\mathcal{V}=
\prod_{n=1}^{\infty}(K_n,\frac{g_n\lambda_n}{Z_n}).$$Then we have:
\begin{corollary}Almost surely in $\mathcal{V}$, $$\lim_{n\to\infty}g_n(x^{(n)})^{\frac{1}{n^2}}=\gamma.$$\end{corollary}
\begin{proof} For all sequences $\{x^{(n)}\}$, $\limsup_ng_n(x^{(n)})^{\frac{1}{n^2}}\leq\gamma$, by hypothesis.If for some $\eta> 0$, $ \liminf_ng_n(x^{(n)})^{\frac{1}{n^2}}\leq \gamma -\eta$, then $x^{(n)} \in B_{\frac{\eta}{2},n}$ for infinitely many $n$ .  Then, by theorem 4.1, and the Borel-Cantelli lemma, the result follows.\end{proof}
Now let $\vec{n}(d)$ be a sequence of multiindices satisfying (3.2).
By corollary 3.6 and theorem 3.9, the hypothesis of theorem 4.1 are satisfied (with $g_n$ replaced by $A_{\vec{n}(d)}$ , $\lambda_n$ by $\mathrm{d}\tau (X)$ and 
$K_n$ by $\Gamma_1^{n_1}\times\ldots\times\Gamma_p^{n_p}\subset L_n$). Corollary 4.2 applies to the product probability space$$\prod_{d=1}^{\infty}(L_{n(d)},\frac{A_{\vec{n}(d)}}{Z_{\vec{n}(d)}}\mathrm{d}\tau (X)).$$By considering the related space$$\prod_{d=1}^{\infty}(\RR^{n(d)},Prob_{\vec{n}(d)})$$
together with (3.6) we have:\begin{corollary}Almost surely in an Angelesco ensemble, a sequence of points $X_{\vec{n}(d)}$ is asymptotically Fekete.\end{corollary}\begin{theorem}Almost surely in an Angelesco ensemble, the  normalized counting measure of a random point converges to the equilibrium measure. That is 
$$\lim_{d\to\infty}\frac{1}{n(d)}\sum\delta (X_{\vec{n}(d)})=\mu_{\Gamma}$$  weak*, almost surely.\end{theorem}\begin{proof}$X_{\vec{n}(d)}$ is almost surely asymptotically Fekete. Hence by corollary 3.4, the result follows.\end{proof} 
\section{the weighted case}
We consider Angelesco ensembles with weights.\\
 That is, in addition to the data for an Angelesco ensemble, namely\\ - $p$ disjoint compact intervals $\Gamma_1,\ldots,\Gamma_p\subset\RR$ \\(We assume that $\Gamma_{j+1}$ lies to the right of $\Gamma_j$ for $j=(1,\ldots,(p-1)$.) \\ -positive Borel measures $\tau_i$ on $\Gamma_i$ for $ i=1,\ldots, p$ \\
-positive reals $r_1,\ldots,r_p$ satisfying $\sum_{i=1}^pr_i=1$.\\-We also have a continuous  real-valued function $Q_i$ on $\Gamma_i$ for $i=1,\ldots,p.$\\We assume that the measures $\tau_i$ satisy the  weighted BM inequality for $Q_i$ on $\Gamma_i$ for $ i=1,\ldots, p$ (see definition 5.3).
The functions $\exp (-Q_i (x))$, continuous positive
functions on $\Gamma_i$,  are the "weights".

Consider the probability measure, invariant under permutation of coordinates, defined analogously to (2.2),(2.3),(2.4) but where\begin{equation}A^Q_{\vec{n}}(X)=\prod_i^p\Delta (X^{(i)})^2\prod_{1\leq i<j\leq n}\Delta(X^{(j)},X^{(i)})\exp(-2n(\sum_{i=1}^p\sum_{j=1}^{n_i}Q_i(x_j^{(i)}))).\end{equation}Since
the probability measure is invariant under permutation of coordinates,it is characterized by\begin{equation}\rho_*(Prob^Q_{\vec{n}})=\frac{A^Q_{\vec{n}}(X)}{Z^Q_{\vec{n}}}d\tau(X)\end{equation} where \begin{equation}Z^Q_{\vec{n}}=\int_{L_n}A^Q_{\vec{n}}(X)d\tau(X).\end{equation}and  $d\tau(X)$ is given by (3.7).\\
The associated Energy functional is  defined as follows: \begin{equation}E^Q(\nu)=E^Q(\nu_1,\ldots,\nu_p)\end{equation}\\$$=\sum_{i=1}^p \int_{\Gamma_i}\int_{\Gamma_i}\log \frac{1}{|x-y|}\mathrm{d}\nu_i(x)\mathrm{d}\nu_i(y)+\sum_{1\leq i<j\leq p}\int_{\Gamma_i}\int_{\Gamma_j}\log\frac{1}{|x-y|}\mathrm{d}\nu_i(x)\mathrm{d}\nu_j(y).$$
$$+2\sum_{i=1}^p\int_{\Gamma_i}Q_i(x)\nu_i(x).$$ The Gonchar-Rachmanov theorem\cite{GR} holds for $E^Q$. That is there is a unique $\mu_{\Gamma,Q}\in\mathcal{M}_r(\Gamma)$ which minimizes $E^Q$ i.e. $E^Q(\mu_{\Gamma,Q})=\inf_{\nu\in\mathcal{M}_r(\Gamma)}E^Q(\nu).$ ( $E^Q$ and $\mu_{\Gamma,Q}$ depend on $\Gamma_1,\ldots,\Gamma_p;r_1,\ldots,r_p;Q_1,\ldots,Q_p$ but not $\tau_1,\ldots,\tau_p.$)  \\ 

Fekete points are again defined by $A^Q_{\vec{n}}\rho (F)\geq A^Q_{\vec{n}}(X)$ for all $X=(X^{(1)},\ldots,X^{(p)})$ with $X^{(i)}\subset\Gamma_i$ for $i=1,\ldots,p.$ Theorem 3.2, lemmas 3.3, 3.4 and corollaries 3.5 and 3.6 are valid in the weighted case with essentially the same proof. Asymptotically Fekete sequences are defined and corollary 3.7 holds. To establish the asymptotic behaviour of the normalizing constants $Z^Q_{\vec{n}(d)}$ in the weighted case, we will need the notion of weighted polynomials and weighted BM inequality.
\begin{definition}Given a weight $\exp(-Q (z))$,where $Q$ is a contiuous real-valued function on a compact set $K\subset \CC$ a \emph{ weighted polynomial} of degree $n$ is a function of the form  $p(z)\exp(-nQ (z))$ where $p$ is an analytic polynomial of degree $\leq n.$\end {definition}
\begin{definition}A positive Borel measure $\tau$ satisfies the\emph{ weighted BM inequality}, for weight  $\exp(-Q (z))$
on a set $K\subset \CC$,  if for all $\epsilon >0$ there is a constant $C>0$ such that $$\sup_{K}|p(z)\exp(-nQ (z))|
\leq C(1+\epsilon )^{deg(p)}\int_K|p(z)|\exp (-nQ (z))\mathrm{d}\tau.$$\end{definition}
\begin{definition}A positive Borel measure satisfies the\emph{ strong weighted BM inequality} if it satisfies the weighted BM inequality on $K$ for all weights $\exp(-Q(z))$ with $Q$ continuous\end{definition}
It is known that for a compact $K\subset\RR$ if $\tau$ satisfies the BM inequality on $K$, then $\tau$ satisfies the strong weighted BM inequality on $K$ (\cite{BB} Theorem 3.2). In particular measures of the form $w(x)\mathrm{d}x$ with $w>0$ a.e. on an interval
satisfy the strong BM inequality on that interval.\\

Now in the weighted case $A^Q_{\vec{n}}$
is, in each variable,  a weighted polynomial of degree $\leq 2n.$ The weighted version of theorem 3.9 follows:
\begin{theorem}$$\lim_d\frac{1}{(n(d))^2}\log Z^Q_{\vec{n}(d)}=-E^Q(\mu_{\Gamma,Q})$$\end{theorem}
 We then obtain the weighted versions of corollary 4.3 and theorem 4.4, specifically
\begin{corollary}Almost surely in a (weighted) Angelesco ensemble, a sequence of points $X_{\vec{n}(d)}$ is asymptotically
Fekete.\end{corollary}
\begin{theorem}Almost surely in a (weighted) Angelesco ensemble,the normalized counting measure of a random point converges to the equilibrium measure.That is $$\lim_{d\to\infty}\frac{1}{n(d)}\delta(X_{\vec{n}(d)})=\mu_{\Gamma,Q}$$ weak*, almost surely.
\end{theorem}

\section{large deviation:statement}
In this section we will state of a large deviation principle (l.d.p.) for Angelesco ensembles. The proof will be given in the next section.\\

Consider a weighted Angelesco ensemble but where each measure $\tau_i$ satisfies the strong weighted BM inequality on $\Gamma_i.$
For $\vec{n}\in\mathbb{N}^p$ we consider the probability measure $Prob_{\vec{n}}^Q$ on $\RR^n$  given by (5.2) and (5.3).Assume that we have a sequence of multiindices $\vec{n}(d)$ satisfying (3.2).

Define $$j^r_{\vec{n}}:\Gamma_1^{n_1}\times\ldots\times\Gamma_p^{n_p}\subset L_n\to \mathcal{M}_r(\Gamma)$$via$$j^r_{\vec{n}}(X)=(\frac{r_1}{n_1}\sum_{j=1}^{n_1}\delta(x_j^{(1)}),\ldots,\frac{r_p}{n_p}\sum_{j=1}^{n_p}\delta(x_j^{(p)})$$
Then, with $\rho$ given by (2.3), $$j^r_{\vec{n}}\circ\rho :\RR^n\to\mathcal{M}_r(\Gamma)$$
We consider measures on $\mathcal{M}_r(\Gamma)$ given by $$\sigma^Q_{\vec{n}}:=(j^r_{\vec{n}}\circ\rho)_*(Prob^Q_{\vec{n}})=(j^r_{\vec{n}})_*\bigg(\frac{A^Q_{\vec{n}}(X)}{Z^Q_{\vec{n}}} d\tau (X)\bigg)$$By definition of the push-forward of a measure this means that, for S a Borel subset of $\mathcal{M}_r(\Gamma)$ then,$$\sigma^Q_{\vec{n}}(S)=\int_{\tilde{S}_{\vec{n}}}
Prob^Q_{\vec{n}}(X),$$\\
where we define$$\tilde{S}_{\vec{n}}=\{X\in L_n|j^r_{\vec{n}}(X)\in S\}.$$

We now define some functionals on $\mathcal{M}_r(\Gamma)$:\\For $G$  open $\subset \mathcal{M}_r(\Gamma)$ we set 
\begin{equation}W^Q_{\vec{n}}(G)=\sup\{ A^Q_{\vec{n}}(X)^{\frac{1}{n^2}}|j^r_{\vec{n}}(X)\in G\}. \end{equation}
\\Now, for the sequence of indices $\vec{n}(d)$ we set\begin{equation}\overline{W}^Q(G)=\limsup_{d\to\infty}W^Q_{\vec{n}(d)}(G);\quad \underline{W}^Q(G)=\liminf_{d\to\infty}W^Q_{\vec{n}(d)}(G)\end{equation} 
\begin{equation}\overline{W}^Q(\mu)=\inf_{G\ni \mu}\overline{W}^Q(G);\quad \underline{W}^Q(\mu)=\inf_{G\ni \mu}\underline{W}^Q(G)\end{equation}.\\

Analagously we set \begin{equation}J^Q_{\vec{n}}(G)=\bigg(\int_{\tilde{G}_{\vec{n}}} A^Q_{\vec{n}}(X)d\tau (X)\bigg)^{\frac{1}{n^2}}\end{equation}
and\begin{equation}\overline{J}^Q(G)=\limsup_{d\to\infty}J^Q_{\vec{n}(d)}(G)  ;\quad\underline{J}^Q(G)=\liminf_{d\to\infty}J^Q_{\vec{n}(d)}(G)\end{equation} 
\begin{equation}\overline{J}^Q(\mu)=\inf_{G\ni \mu}\overline{J}^Q(G);\quad \underline{J}^Q(\mu)=\inf_{G\ni \mu}\underline{J}^Q(G)\end{equation}.\\
For $Q=(0,\ldots,0)$ we use the notation $ \overline{W}(\mu),\overline{J}(\mu)$ etc...\\
\begin{proposition}$\overline{J}^Q(\mu)\leq\overline{W}^Q(\mu);\quad\underline{J}^Q(\mu)\leq\underline{W}^Q(\mu);\\ \quad\log\overline{W}^Q(\mu)\leq-E^Q(\mu)$\end{proposition}
\begin{proof}The first two inequalities follow from the fact that$$J_{\vec{n}}^Q(G)\leq W^Q_{\vec{n}}(G)$$
the last inequality follows from the fact that $-E^Q$ is u.s.c.\end{proof}
\begin{definition}
A sequence $\{\beta_d\}$ of probability measures on a compact Hausdorff space $\mathcal{X}$ satisfies a l.d.p. with good rate function $R$ and speed ${a_d}$ if:\\(1) ${a_d}$ is a sequence of positive reals $\{a_d\}\to 0$.\\(2) $ R$ is a non-negative real-valued lower semicontinuous function on $\mathcal{X}$\\ (3) for all Borel sets $\mathcal{B}\subset\mathcal{X}$we have$$-\inf_{x\in  \mathcal{B}^{\circ}}R(x)\leq\liminf_{d\to\infty}\frac{1}{a_d}\log \beta_d(\mathcal{B})\leq\limsup_{d\to\infty}\frac{1}{a_d}\log \beta_d(\mathcal{B})\leq-\inf_{x\in\overline{\mathcal{B}}} R(x).$$\end{definition}
	It is known(see \cite{DZ},theorem 4.1.11)that for $\mathcal{G}$ a base for the topology of $\mathcal{X}$ that there is an l.d.p. with good rate function R if \begin{equation}\inf_{\mathcal{G}(x)}(\liminf_{d\to\infty}\frac{1}{a_d}\log \beta_d(G))=\inf_{\mathcal{G}(x)}(\limsup_{d\to\infty}\frac{1}{a_d}\log \beta_d(G))\end{equation} where $\mathcal{G}(x)$ denotes the collection of sets in $\mathcal{G}$ which contain $x\in\mathcal{X}$\\$-R(x)$ is then the common value of the two sides of (6.7)

\begin{theorem}Let $\vec{n}(d)$ be a sequence satisfying (3.2).The sequence of measures $ \sigma^Q_{\vec{n}(d)}$ on $\mathcal{M}_r(\Gamma)$ satisfies a large deviation principle with speed $(n(d))^2$and good rate function$$R(\mu)=\log W^Q(\mu_{\Gamma,Q})-\log W^Q(\mu)=E^Q(\mu)-E^Q(\mu_{\Gamma,Q}).$$\end{theorem}
Now, by definition,$$\frac{1}{n^2}\log \sigma^Q_{\vec{n}}(G)=\log J^Q_{\vec{n}}(G)-\frac{1}{n^2}\log Z^Q_{\vec{n}}.$$
so, using (6.7) and theorem 5.4, to prove theorem 6.3, it suffices to show that for all measures $\mu\in\mathcal{M}_r(\Gamma)$\begin{equation}\log \overline{J}^Q(\mu)=\log\underline{J}^Q(\mu)=\log \overline{W}^Q(\mu)=\log\underline{W}^Q(\mu)=-E^Q(\mu).\end{equation}
This will be done in the next section.
\section{large deviation:proofs}
\begin{proposition}$$\overline{W}(\mu)=\overline{W}^Q(\mu)\exp (2\sum_{i=1}^p\int_{\Gamma_i}Q_id\mu_i)$$\end{proposition}
\begin{proof}Given $\epsilon>0$,there is  a neighbourhood, $G$ of $\mu$ in $\mathcal{M}_r(\Gamma)$ such that$\sum_{i=1}^p|\int_{\Gamma_i}Q_i(d\mu_i-d\theta_i)| \leq\epsilon$ for all $\theta =(\theta_1,\ldots,\theta_p)\in G$. Thus,$$-\epsilon\leq\sum_{i=1}^p\int_{\Gamma_i}Q_id\mu_i-\sum_{i=1}^p\frac{r_k}{n_k}\sum_{k=1}^{n_i}\delta (x_k^{(i)}) \leq\epsilon$$ for all $X=(X^{(1)},\ldots,X^{(p)})\in\tilde{G}_{\vec{n}}.$\\It follows  that for all such $X$ and all $n(d)$,sufficiently large, we have $$-2\epsilon\leq\sum_{i=1}^p\int_{\Gamma_i}Q_id\mu_i-\frac{1}{n}\sum_{i=1}^p\sum_{k=1}^{n_i}Q_i(x_k^{(i)})\leq2\epsilon.$$
This implies that for $X\in\tilde{G}_{\vec{n}}$, we have$$e^{-4n^2\epsilon}A_{\vec{n}}^Q(X)\leq \exp (2n^2(\sum_{i=1}^p\int_{\Gamma_i}Q_id\mu_i))A_{\vec{n}}(X)\leq e^{4n^2\epsilon}A_{\vec{n}}^Q(X).$$
Taking the $\frac{1}{n^2}$ power and then sup over $X\in\tilde{G}_{\vec{n}}$ we obtain$$e^{-2\epsilon}W_{\vec{n}}^Q(G)\leq \exp(2(\sum_{i=1}^p\int_{\Gamma_i}Q_id\mu_i))W_{\vec{n}}(G)\leq e^{2\epsilon}W_{\vec{n}}^Q(G).$$
Then taking the $\limsup$ as $d\to\infty$ we obtain$$e^{-2\epsilon}\overline{W}^Q(G)\leq \exp(2(\sum_{i=1}^p\int_{\Gamma_i}Q_id\mu_i))W(G)\leq e^{2\epsilon}\overline{W}^Q(G).$$
Now, for any $\epsilon>0$, a neighbourhood of $\mu$ may be found so the above inequality holds. Hence, taking the $\inf$ over all neighbourhoods of $\mu$ the result follows.
\end{proof} 
\begin{corollary}
$$\underline{W}(\mu)=\underline{W}^Q(\mu)\exp (2\sum_{i=1}^p\int_{\Gamma_i}Q_id\mu_i).$$
$$\overline{J}(\mu)=\overline{J}^Q(\mu)\exp (2\sum_{i=1}^p\int_{\Gamma_i}Q_id\mu_i).$$
$$\underline{J}(\mu)=\underline{J}^Q(\mu)\exp (2\sum_{i=1}^p\int_{\Gamma_i}Q_id\mu_i).$$
\end{corollary}
Let $G$ be a neighbourhood of $\mu_{\Gamma,Q}\in\mathcal{M}_r(\Gamma)$ and let $\delta ^Q$ be defined by \begin{equation}\log\delta ^Q=-E^Q(\mu_{\Gamma,Q}).\end{equation}\\
For $\eta>0$, let$$B^Q_{\eta,\vec{n}}:=\{X\in \Gamma_1^{n_1}\times\ldots\times\Gamma_p^{n_p}|A^Q_{\vec{n}}(X)\leq(\delta^Q-\eta)^{n^2}\}$$ and$$(B^Q_{\eta,\vec{n}})^c:=\{X\in \Gamma_1^{n_1}\times\ldots\times\Gamma_p^{n_p}|A^Q_{\vec{n}}(X)>(\delta^Q-\eta)^{n^2}\}$$ Note that this notation is compatible with that used in theorem 4.1.
\begin{proposition}Let $\vec{n}(d)$ be a sequence of multiindices satisfying (3.2). Given a sequence $\eta_m\downarrow  0$ and any sequence of integers $d_m\uparrow\infty$ there exists a $m_0$ such that for $m\geq m_0$ we have$$(B^Q_{\eta_m,\vec{n}(d_m)})^c\subset\tilde{G}_{\vec{n}(d_m)}. $$\end{proposition}  
\begin{proof}The proof will be by contradiction.\\Suppose not, then for an infinite sequence of $m's$ we have a point$$
{}_mX\in ( B^Q_{\eta_m,\vec{n}(d_m)})^c\setminus\tilde{G}_{\vec{n}(d_m)}.$$Since $\eta_m\downarrow 0,$   we have ${}_mX$ is asymptotically Fekete so $\lim_m \delta(\rho({}_mX))=\lim_mj^r_{\vec{n}}(d_m)(X)=\mu_{\Gamma,Q}$ weak*. This contradicts ${}_mX\not\in G$ for that sequence of $m's$.
\end{proof}
\begin{theorem}$$\log\underline{J}^Q(\mu_{\Gamma,Q})=\log\delta^Q.$$\end{theorem}  
\begin{proof}We may assume that $\delta^Q>0$, since each of the intervals $\Gamma_i$  is non-degenerate\cite{BKMW})\\Given a sequence $\eta_m\downarrow 0$ we can choose $d_m$ so that, using theorem 4.1 we have $$\rho_*(Prob_{\vec{n}(d_m)})(B^Q_{\eta_m,\vec{n}(d_m)})\leq \bigg (1-\frac{\eta_m}{2\delta^Q}\bigg )^{(n(d_m))^2)}$$Passing to a subsequence of the $d_m's$ we may assume that$$\lim_m\bigg  (1-\frac{\eta_m}{2\delta^Q}\bigg  )^{(n(d_m))^2}= 0$$Now, for$m\geq m_0$ (given by proposition 7.3), $$J^Q_{\vec{n}(d_m)}(G)=\int_{\tilde{G}_{\vec{n}(d_m)}}A^Q_{\vec{n}(d_m)}(X)d\tau (X)\geq\int_{(B_{\eta_m,\vec{n}(d_m)})^c}A^Q_{\vec{n}(d_m)}(X)d\tau (X)$$which is$$\geq   {Z^Q_{\vec{n}(d_m)}}\bigg(1-\bigg(1-\frac{\eta_m}{2\delta^Q}\bigg)^{(n(d_m))^2}\bigg)$$Now using the asymptotics for $Z^Q_{\vec{n}}$ (i.e. theorem 5.3) we conclude that$$\liminf_m\log\underline{J}^Q_{\vec{n}(d_m)}\geq\log \delta^Q$$
 This relation, in fact, holds taking $\liminf$ over the full sequence of integers.Thus $$\log\underline{J}^Q(\mu_{\Gamma,Q})=\log \delta^Q$$\end{proof}
\begin{corollary}Let $H=(H_1,\ldots,H_p)$ where $H_i$ is continuous on $\Gamma_i.$Then (6.8) holds with $\mu=\mu_{\Gamma,H}$\end{corollary}
\begin{proof}By corollary 7.2, $\log\underline{J}^Q(\mu)-\log\underline{J}^H(\mu)=2\sum_{i=1}^p\int_{\Gamma_i}(H_i-Q_i)d\mu=
E^H(\mu)-E^Q(\mu)$
\end{proof}
It remains to show that (6.8) holds for a general $\mu \in \mathcal{M}_r(\Gamma)$. If $E^Q(\mu)=+\infty$ then the inequalities of Proposition (6.1) show that (6.8) holds. We may therefore assume that $E^Q(\mu)<+\infty.$ The proof will be completed by showing that an arbitrary such measure can be approximated by weighted equilibrium measures in an appropriate way. Specifically we have:\begin{theorem}Let $\mu \in \mathcal{M}_r(\Gamma)$ and suppose that $E^Q(\mu) < +\infty.$ Then there exists a sequence of weighted equilibrium measures $\mu_{\Gamma,H_m}$
with $\mu_{\Gamma,H_m}\rightarrow \mu$ weak* in $\mathcal{M}(\Gamma)$ and $E^Q(\mu_{\Gamma,H_m})\rightarrow E^Q(\mu).$\end{theorem}
\begin{proof}(of (6.8) from theorem 7.6 )\\
Given $\mu$ we take a sequence of equilibrium measures (as in theorem 7.6 ) $\mu_{\Gamma,H_m}$.Then using the u.s.c. property of the functional, we have $$\log \underline{J}^Q(\mu)\geq \limsup_m\log\underline{J}^Q(\mu_{\Gamma,H_m})=-\limsup_mE^Q(\mu_{\Gamma,H_m})=-E^Q(\mu).$$The middle inequality is due to theorem 7.4. The reverse inequality having been noted in proposition 6.1, (6.8) follows.\end{proof}
\begin{proof}(of theorem 7.6) For $\alpha,\beta$ two measures with compact support in the complex plane we consider their mutual energy:$$I(\alpha,\beta)=-\int\int\log |x-y|d\alpha d\beta.$$We also let $$p_{\alpha}(y)=-\int\log |x-y|d\alpha (x)$$\\Now for $\mu=(\mu_1,\ldots,\mu_p)\in\mathcal{M}(\Gamma)$ we use the partial potentials, defined for $s=1,\ldots,p$ by$$U^{\mu}_s=\frac{1}{2}\sum_{j=1}^pp_{\mu_j}+\frac{1}{2}p_{\mu_s}.$$We have $$E^Q(\mu)=\sum_{i=1}^pI(\mu_i,\mu_i)+\sum_{1\leq i<j\leq p}I(\mu_i,\mu_j)+2\sum_{i=1}^p\int  Q_id\mu_i.$$Then direct calculation gives $$E^Q(\nu )-E^Q(\mu )=2\sum_{i=1}^p\int  (U^{\mu}_i+Q_i)d(\nu_i -\mu_i)+E(\nu -\mu)$$ where for $Q=(0,\ldots,0)$ we use the notation $E$. This shows that if $U_i^{\mu}$ continuous and $Q_i=-U_i^{\mu}$ for $i=1,\ldots,p$, then the unique minimum of the left side (as a function of $\nu$) is attained when $\nu=\mu$.\\

Now consider $ -p_{\mu_j}$ which is usc on $\Gamma_j$ for $j=1,\ldots,p$. Note that since $E^Q(\mu)<+\infty$ that $I(\mu_i,\mu_i)<+\infty$ for $i=1,\ldots,p.$ By(\cite{BL} Theorem 6.2), there exists a sequence of measures $\mu_{m,j}$ for each $j=1,\ldots,p$ converging to $\mu_j$
weak* and with $ I(\mu_{m,j},\mu_{m,j})$ converging to $ I(\mu_j,\mu_j).$  It follows that $ I(\mu_{m,j},\mu_{m,i})$ converges to $I(\mu_j,\mu_i)$ for all $i,j$, since the convergence for $i\neq j$ follows from the weak convergence of the measures and the fact that $\Gamma_i$ and $\Gamma_j$ are disjoint. Hence $E^Q(\mu_m)$ converges to $E^Q(\mu)$ for any $Q$. The sequence $\mu_{m,j}$ can be chosen (see \cite{BL}) so that the potentials $ p_{\mu_{m,j}}$ are continuous in the plane and so the partial potentials are continuous. Thus each $\mu_m$ is an equilibrium measure.
This completes the proof of theorem 7.6.
\begin{remark}Weak* convergence of a sequence of measures $\mu_m$ to $\mu$ does not, in general, imply the convergence of $I(\mu_m,\mu_m)$ to $I(\mu,\mu)$. For an example see \cite{BG}.\end{remark}

\end{proof}

\end{document}